\documentclass[12pt,titlepage]{article}
\usepackage{graphicx,color}
\usepackage{amsmath,amssymb,amsthm}
\input epsf

\newtheorem{thm}{Theorem}

\newtheorem{cor}{Corollary}

\newtheorem{prop}{Proposition}

\newtheorem{conj}{Conjecture}
\newtheorem{defi}{Definition}

\newtheorem{remark}{Remark}

\def\beq{\begin{equation}}\def\eeq{\end{equation}}
\def\beqn{\begin{eqnarray}}\def\eeqn{\end{eqnarray}}

\def\qed{\ifhmode\unskip\nobreak\fi\quad\ifmmode\Box\else$\Box$\fi}

\title{Partitioning transitive tournaments into isomorphic digraphs}

\author{{\bf Attila Sali}$^{a,b,}$\thanks{Research partially supported by the
    National Research, Development and Innovation Office (NKFIH)
    grant K--116769. This work is also connected to the scientific program of the "Development of quality-oriented and harmonized R+D+I strategy and functional model at BME" project, supported by the New Hungary Development Plan (Project ID: TÁMOP-4.2.1/B-09/1/KMR-2010-0002).
}
$\qquad$ {\bf G\'abor Simonyi}$^{a,b,}$\thanks{Research partially supported by the
    National Research, Development and Innovation Office (NKFIH)
    grants K--116769 and K--120706. This work is also connected to the scientific program of the "Development of quality-oriented and harmonized R+D+I strategy and functional model at BME" project, supported by the New Hungary Development Plan (Project ID: TÁMOP-4.2.1/B-09/1/KMR-2010-0002).
}
$\qquad$ {\bf G\'abor Tardos}$^{a,}$\thanks{Research partially supported by the Cryptography ``Lend\"ulet'' project of the Hungarian Academy of Sciences and by the National Research, Development and Innovation Office (NKFIH) grants K--116769 and SSN-117879.}\\ \\
$^a$Alfr\'ed R\'enyi Institute of Mathematics, Hungarian Academy of
Sciences\\ \\
$^b$Department of Computer Science and Information Theory,\\ Budapest University of Technology and Economics\\ \\
 {\tt sali@renyi.hu} \ \ \  {\tt simonyi@renyi.hu} \ \ \ {\tt tardos@renyi.hu}}
\date{}

\begin{document}
\maketitle
\begin{abstract}
In an earlier paper (see \cite{SaSi}) the first two authors have shown that
self-complementary
graphs can always be oriented in such a way that the union of the oriented
version and its isomorphically oriented complement gives a transitive
tournament. We investigate the possibilities of generalizing this theorem to
decompositions of the complete graph into three or more isomorphic graphs.
We find that a complete characterization of when an orientation with similar
properties is possible seems
elusive. Nevertheless, we give sufficient conditions that generalize the earlier theorem and
also imply that decompositions of odd vertex complete graphs to Hamiltonian
cycles admit such an orientation. These conditions are further generalized and
some necessary conditions are given as well.

\end{abstract}

\section{Introduction}
\message{Introduction}

Investigating the relationship between the Shannon capacity of  graphs and the Sperner capacity of their oriented versions the first two authors proved the following theorem.

\begin{thm} \label{thm:selfc} {\rm (\cite{SaSi}, cf. also \cite{GyA})} Let $G$ be a graph isomorphic to its complement $F=\overline{G}$. Then $G$ and $F$ can be oriented so that they remain isomorphic as digraphs while the tournament formed by their union is the transitive tournament.

Moreover, the above can be done for any fixed isomorphism between $G$ and $F$. That is, for any such isomorphism $f$ one can find oriented versions $\vec G$ and $\vec F$ of $G$ and $F$, respectively, such that $f$ provides and isomorphism between $\vec G$ and $\vec F$ and the union of $\vec G$ and $\vec F$ is a transitive tournament.
\end{thm}

The goal of the present paper is to investigate the possibilities of
generalizing the above theorem to three or more graphs, that is, to the
situation when (the edge set of) the complete graph is partitioned into three
or more isomorphic graphs.
As already observed by G\"orlich, Kalinowski, Meszka, Pil\'sniak, and Wo\'zniak \cite{GKMPW} in this case it will not be true that the three graphs can always be oriented
in an isomorphic manner so that their union forms a transitive
tournament. Moreover, a complete characterization of when this is possible
seems to be elusive. In \cite{GKMPW, GKMPW2} the authors determine all digraphs with at most four edges that can decompose a transitive tournament. (For related results, see also \cite{GPW} and \cite{GP}.)

In our approach we fix the number of isomorphic graphs in the decompositions considered. We start with the case when this number is $3$. We will give some sufficient conditions when the isomorphic parts of a decomposition of $K_n$ can be isomorphically oriented to get a decomposition of the transitive tournament. This result gives a generalization of Theorem~\ref{thm:selfc}.

It is well-known that complete graphs on an odd number of vertices decompose
into Hamiltonian cycles. One can directly show how to obtain a decomposition of the transitive tournament of odd order into isomorphically oriented Hamiltonian cycles, but this also follows from our sufficient condition mentioned above.

First we extend our sufficient condition to a more general one, then we show with an example (found by computer) that this more general condition is still not necessary. A
complete characterization seems out of reach, but we are able to give some non-trivial necessary conditions.


As usual $K_n$ denotes the complete graph on $n$ vertices, while we denote the transitive tournament on $n$ vertices by $T_n$. The vertex set of $K_n$ and $T_n$ is assumed to be $[n]=\{0,1\dots,n-1\}$ and we consider these vertices as residue classes modulo $n$, that is equality between vertices will be understood modulo $n$. We denote the cyclic permutation of $[n]$ bringing $i$ to $i+1$ by $\sigma_n$.

\section{Small examples and problem formulation}\label{sect:smallex}

First we recall an example from \cite{GKMPW} of three isomorphic graphs partitioning the complete graph that cannot be isomorphically oriented so that their union is a transitive tournament even if the functions giving the isomorphism among them are not fixed.

Let $n=|V(G)|=4$ and the three isomorphic graphs be paths on $3$ vertices. One
easily partitions $K_4$ into three such graphs. It is also easy to see that
whatever way we orient these paths in an isomorphic manner, we cannot put them
together to obtain a transitive tournament on $4$ vertices. This is simply because from no orientation can we produce simultaneously a vertex of outdegree $0$ and a vertex of outdegree $3$. Note that this example is just a very special case of Theorem 5 from \cite{GKMPW}.

Let us assume that the edge set of the complete graph $K_n$ is partitioned into three isomorphic graphs $F$, $G$ and $H$. We can ask whether there are isomorphic orientations $\vec F$, $\vec G$ and $\vec H$ of the graphs $F$, $G$ and $H$, respectively, such that their union gives a transitive tournament. But we can be more specific and fix an isomorphism $\sigma$ from $F$ to $G$, an isomorphism $\rho$ from $G$ to $H$ and ask whether there are orientations $\vec F$, $\vec G$ and $\vec H$ of the graphs $F$, $G$ and $H$ whose union is a transitive tournament and such that $\sigma$ is an isomorphism between $\vec F$ and $\vec G$ and $\rho$ is an isomorphism between $\vec G$ and $\vec H$.

To illustrate the difference between these two questions let us consider the smallest possible example. The graph $K_3$ can be partitioned into three (isomorphic) single edge graphs: $F$, $G$ and $H$. Clearly, the three oriented edges of $T_3$ also form
isomorphic graphs. This answers the first question above for this specific partition affirmatively. If, however, we fix a cyclic permutation $\sigma=\rho$ that brings $F$ to $G$ and $G$ to $H$, then the answer to the second question is negative. Indeed, if $\vec F$ is either orientation of $F$, $\vec G=\sigma(\vec F)$ and $\vec H=\sigma^2(\vec F)$, then the union of these three directed graphs is a directed cycle and thus not transitive.

In this paper we will concentrate on the question with fixed permutations $\sigma$ and $\rho$. We will only consider the special case $\sigma=\rho$. Although this assumption is restrictive, it is
in complete analogy with the case of two self-complementary graphs, and we
believe that understanding this special case would largely improve our
knowledge about the situation.

\begin{defi}
Let $\sigma$ be a permutation of the vertex set of $K_n$. We call the partition of the edge set of $K_n$ into three graphs $F$, $G$ and $H$ such that the permutation $\sigma$ brings $F$ to $G$ and $G$ to $H$ a {\em $\sigma$-partition}. In this case $\sigma$ brings $H$ to $F$ and $\sigma^3$ is an automorphism of all three of these graphs. We call a transitive orientation $T$ of $K_n$ a {\em transitive $\sigma$-orientation} of this $\sigma$-partition if the subgraphs $\vec F$, $\vec G$ and $\vec H$ of $T$ that are the orientations of the graphs $F$, $G$ and $H$, respectively, satisfy $\sigma(\vec F)=\vec G$ and $\sigma(\vec G)=\vec H$. We say that $\sigma$ {\em reverses the orientation} of an edge $e$ in $T$ if $\sigma(e)$ is oriented in the other direction, that is, if $e$ goes from $a$ to $b$, and the edge of $T$ between $\sigma(a)$ and $\sigma(b)$ is oriented toward $\sigma(a)$. Observe that a transitive orientation $T$ of $K_n$ is a transitive $\sigma$-orientation of the $\sigma$-partition of $K_n$ to $F$, $G$ and $H$ if and only if $\sigma$ reverses no edges of $T$ that belong to $F$ or $G$.
\end{defi}

Just as it was the case with
self-complementary graphs, we may assume that the permutation $\sigma$ is cyclic as the case of general $\sigma$ reduces to the cyclic case. Indeed, let the cycle decomposition of the permutation of $\sigma$ be $\rho_1\rho_2\dots\rho_k$. Let $F$, $G$ and $H$ form a $\sigma$-partition $P$. The subgraphs $F_i$, $G_i$, $H_i$ induced by the domain of the cycle $\rho_i$ form a $\rho_i$-partition $P_i$ for all $i$. Clearly, if $P$ has a transitive $\sigma$-orientation, then it restricts to transitive $\rho_i$-orientations of $P_i$. On the other hand, if $P_i$ has a transitive $\rho_i$-orientation for all $i$, then $P$ has a transitive $\sigma$-orientation. To see this last statement simply keep the orientations of the edges in the transitive $\rho_i$-orientations and orient edges connecting vertices from distinct cycles $\rho_i$ and $\rho_j$ toward the higher indexed cycle. None of these latter type of edges is reversed by $\sigma$.

It is easy to see that a $\sigma$-partition exists if and only if $\sigma$ has at most one fixed point and the length of all non-trivial cycles of the cycle decomposition of $\sigma$ is divisible by $3$.

From now on we do make the assumption that $\sigma$ consists of a single cycle on $n>1$ vertices,
namely $\sigma=\sigma_n$, where $\sigma_n$ stands for the permutation on the set $[n]=\{0,1,2,\dots,n-1\}$ bringing $i$ to $i+1\bmod n$. The vertices of our graphs will therefore be the elements of $[n]$ and we consider them as the residue classes modulo $n$, that is, equalities about them are always understood modulo $n$. We assume $n$ is divisible by $3$ as otherwise there is no $\sigma_n$-partition.

We denote the graphs of the $\sigma_n$-partition by $F_0$, $F_1=\sigma_n(F_0)$ and $F_2=\sigma_n^2(F_0)$. The {\em label} $\ell(a,b)$ of an edge $\{a,b\}$ of $K_n$ is the index of the subgraph the edge belongs to, so the label of the edges of $F_i$ are $i$. As $\sigma_n$ brings $F_0$ to $F_1$ to $F_2$ and back to $F_0$ we must have $\ell(a+1,b+1)\equiv\ell(a,b)+1$ for all $a$ and $b$, where the congruence is modulo $3$ (and, as noted above, the vertices are understood modulo $n$). With the same convention we have the more general congruence for any edge $\{a,b\}$ and integer $i$:
\begin{equation}\label{shift}
\ell(a+i,b+i)\equiv\ell(a,b)+i\pmod3.
\end{equation}

\begin{defi}
The {\em defining sequence} of the $\sigma_n$-partition $\{F_0,F_1,F_2\}$ is $a_1,a_2,\dots,a_m$, where $m=\lfloor n/2\rfloor$ and $a_j=\ell(0,j)$. By the congruence above, this sequence determines all other labels and thus the entire $\sigma_n$-partition. On the other hand, it is easy to see that (as $n$ is divisible by $3$) every sequence of length $\lfloor n/2\rfloor$ over the alphabet $\{0,1,2\}$ is a defining sequence of a $\sigma_n$-partition. This is
analogous to the case of self-complementary graphs, cf. \cite{english, Ringel, Sachs}. By symmetry, we may and will often assume that $\ell(0,1)=0$, that is, the defining sequence starts with $a_1=0$. This can be achieved by shifting the $\sigma_n$-partition by $\sigma_n$ or $\sigma_n^2$.
\end{defi}

In the smallest $n=3$ case, there is just one $\sigma_3$-partition and we have already seen that it has no transitive $\sigma_3$-orientation. Let us look at the next case $n=6$ a bit closer. By the
foregoing, there are $3^2=9$ $\sigma_6$-partitions to consider according to the labeling of the edges $\{0,2\}$ and $\{0,3\}$.

The corresponding graphs $F_0$ are depicted in Figure~1. It turns out that transitive $\sigma_6$-orientations exist in
exactly four of the nine cases. (We have indicated such an orientation in Figure~1 whenever it exists.)
Notice that the $F_0$ is simply a path on the six vertices in four cases
but a transitive $\sigma_6$-orientation exists for only two of
them. (The truth of this statement will follow from the results of the next
section.) Thus, in spite of the isomorphism of these four graphs they behave
differently according to the different effect of permutation $\sigma_6$ on them.

\begin{figure}[!htbp]
\centering
\includegraphics[scale=0.56]{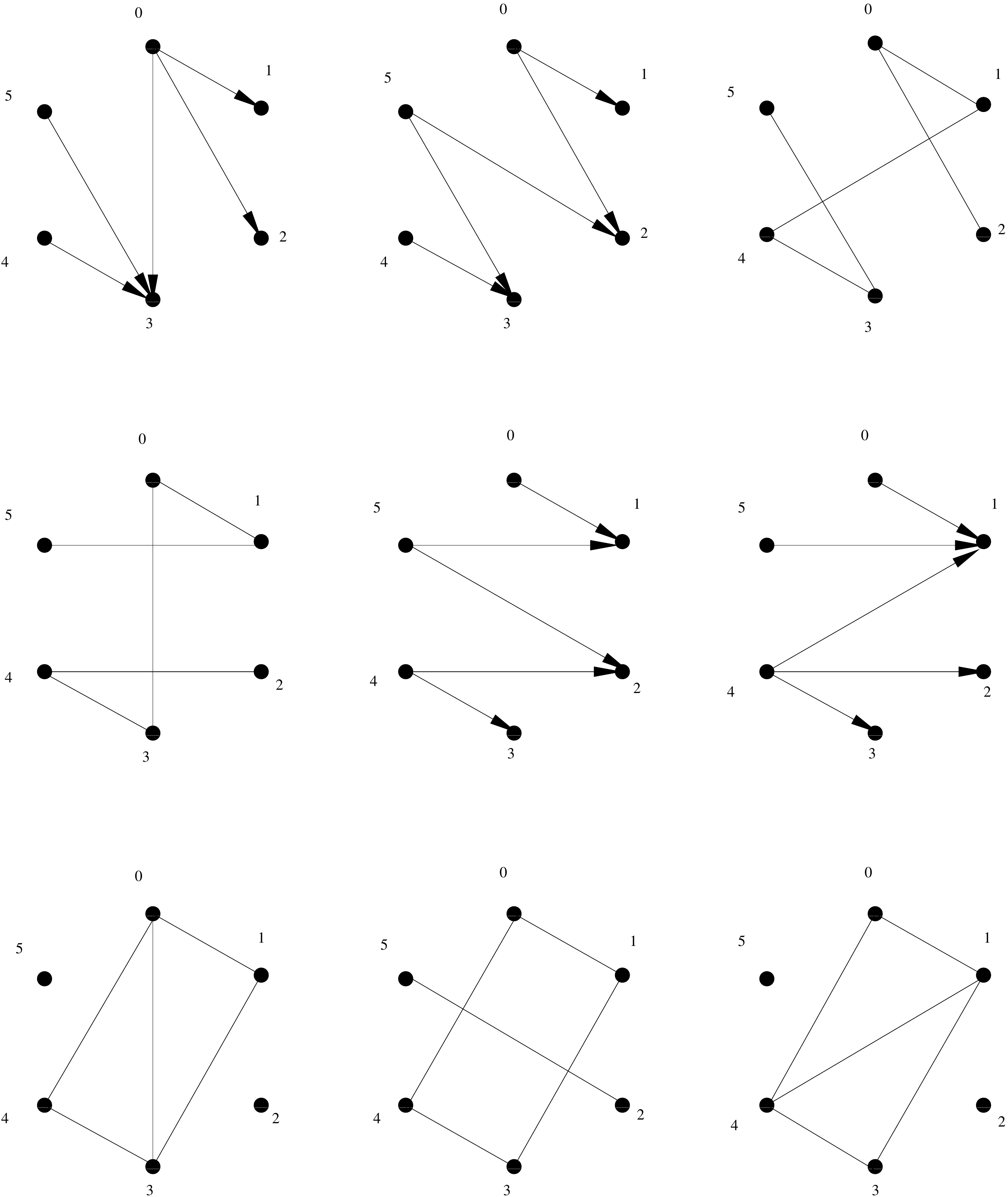}
\caption{Four of the nine possible $3$-partitions of $K_6$ can be oriented as
  required. In the remaining five cases such orientations do not exist.}
\label{fig:kilenc}
\end{figure}

\section{The standard orientation}\label{sect:standard}

We want to decide whether a given $\sigma_n$-partition $P$ has a transitive
$\sigma_n$-orientation. For our notation including the labeling of edges and
the definition of the defining sequence see the previous section. We will
describe a transitive $\sigma_n$-orientation with an ordering
$\tau(1),\tau(2),\dots,\tau(n)$ on the vertices. We say that an orientation is
{\em consistent} with $\tau$ if all edges point towards the vertex that come
later in the order. Clearly, the transitive orientation of $K_n$ and the
ordering it is consistent with mutually determine each other. Recall that a
transitive orientation $T$ is a transitive $\sigma_n$-orientation of our
$\sigma_n$-partition $P$ if and only if the orientation of no label-$0$ or
label-$1$ edge is reversed by $\sigma_n$.

\begin{thm}\label{standard}
If $n=2m$ and the defining sequence $a_1\dots a_m\in\{0,1,2\}^m$ of a
$\sigma_n$-partition satisfies that for every $j\in\{1,\dots,m-1\}$ either
$a_{j+1}=a_j$ or $a_{j+1}\equiv a_j+1 \pmod3$, then there exists a transitive $\sigma_n$-orientation for this $\sigma_n$-partition.
\end{thm}

\proof
The proof is an extension of the argument given by Gy\'arf\'as \cite{GyA} for our Theorem~\ref{thm:selfc}. We give a linear order $\tau$ of the vertices $0,1,\dots,n-1$ and show that orienting the edges consistently with this order gives a transitive $\sigma_n$-orientation.

Let us first recall our assumption that $a_1=0$. This can be achieved by appropriately relabeling the vertices. The relabeling changes the defining sequence but does not affect the condition in the theorem. Now we define $\tau$. We set $\tau(1)=0$ and declare that $\tau$ will have the property, that for any $i$, the set of vertices $A_i:=\{\tau(1), \tau(2),\dots,\tau(i)\}$ forms a consecutive arc of the cycle formed by the vertices $0,1,\dots,n-1$, i.e., it is equal to $\{j_i+1,\dots,j_i+i\}$ for some $j_i\in \{n~-~i,n~-~i~+~1,\dots,n~-~1\}$. Recall that the names of the vertices are understood modulo $n$. Now
$\tau(i+1)$, that is the unique element of $A_{i+1}\setminus A_i$ is either
$j_i$ or $j_i+i+1$ for every $i$.
Thus $\tau$ is determined if we give a rule for deciding which of the two elements $j_i$ and $j_i+i+1$ should be taken as $\tau(i+1)$ if $i<n-1$. (No rule is needed for $i=n-1$ as then $j_i=j_i+i+1$ is the only vertex outside $A_i$.) This choice for $\tau(i+1)$ depends on the label of the edge $\{j_i,j_i+i+1\}$. If $\ell(j_i,j_i+i+1)=0$, then we set $\tau(i+1)=j_i$ making $j_{i+1}=j_i-1$. If $\ell(j_i,j_i+i+1)=2$, then we set $\tau(i+1)=j_i+i+1$ making $j_{i+1}=j_i$. We claim that the third possibility, namely $\ell(j_i,j_i+i+1)=1$ will not happen for any $1\le i<n-1$.

First we show this last statement by induction. Note that all congruences are modulo $3$. The base case is all right as
$\tau(1)=0$ and $\ell(n-1,1)\equiv a_2-1$, see Equation~(\ref{shift}). By the assumption on the defining sequence this is either $a_1-1\equiv2$ or $(a_1+1)-1=0$. Now assume the statement to be true for the edge $\{j_{i-1},j_{i-1}+i\}$, and we show that it is also true for $\{j_i,j_i+i+1\}$. By Equation(\ref{shift}), we have $\ell(j_{i-1},j_{i-1}+i)\equiv\ell(0,i)+j_{i-1}$ and $\ell(j_i,j_i+i+1)\equiv\ell(0,i+1)+j_i$. We have $j_i=j_{i-1}-1$ if $\tau(i)=j_{i-1}$, i.e., if $\ell(j_{i-1},j_{i-1}+i)=0$, while otherwise this label is $2$, so we have $j_i=j_{i-1}$. Therefore, we can formulate the congruence:
$$j_i\equiv j_{i-1}-\ell(j_{i-1},j_{i-1}+i)-1\equiv j_{i-1}-(\ell(0,i)+j_{i-1})-1=-\ell(0,i)-1.$$
We also have:
\begin{equation}\label{bar}
\ell(j_i,j_i+i+1)\equiv\ell(0,i+1)+j_i\equiv\ell(0,i+1)-\ell(0,i)-1.
\end{equation}
For $1\le i\le n/2-1$ we simply have $\ell(0,i+1)=a_{i+1}$ and $\ell(0,i)=a_i$, so
\begin{equation}\label{kicsi}
\ell(j_i,j_i+i+1)\equiv a_{i+1}-a_i-1,
\end{equation}
so by assumption it cannot be $1$. If $n/2\le i<n$ we use Equation~(\ref{shift}) again to see that $\ell(0,i)=\ell(i,0)\equiv\ell(0,n-i)+i=a_{n-i}+i$ and similarly, $\ell(0,i+1)\equiv a_{n-i-1}+i+1$. Therefore \begin{equation}\label{nagy}
\ell(j_i,j_i+i+1)\equiv a_{n-i-1}-a_{n-i},
\end{equation}
which cannot be $1$ either by the same assumption.
Note that we used the fact that $n$ is even. For $n$ odd and $i=(n-1)/2$ we would have $\ell(j_i,j_i+i+1)\equiv\ell(0,i+1)-\ell(0,i)-1\equiv(a_i-i)-a_i-1=-i-1\equiv1$.

\medskip
\par\noindent
We need to show that the orientation consistent with the order $\tau$ is a transitive $\sigma_n$-orientation. As noted above, for this we have to show that $\sigma_n$ reverses the orientation only of edges of label $2$. Equivalently, if an edge $\{u,v\}$ is oriented from $u$ to $v$ and
it has label $1$ or $2$ then the edge $\{u-1,v-1\}$ is oriented from
$u-1$ to $v-1$. Assume $\{u,v\}$ is oriented from $u$ to $v$, that is, $\tau^{-1}(u)<\tau^{-1}(v)$. We distinguish cases according to the order of $u$ and $v$. Note that while in most formulas we consider the vertices as residue classes modulo $n$ (and thus equality really means congruence modulo $n$) in inequalities the vertices are treated as integers between $0$ and $n-1$.

In the simplest case we have $0<u<v$. In this case $u-1\in A_{\tau^{-1}(u)}$ while
$v\notin A_{\tau^{-1}(u)}$ and either $(v-1)\notin A_{\tau^{-1}(u)}$ or
$v-1=u$. In both cases we have $\tau^{-1}(u-1)<\tau^{-1}(v-1)$
implying that the edge $\{u-1,v-1\}$ is oriented from $u-1$ to $v-1$ as we
need.

If $u=0$ and $\{u-1,v-1\}$ is oriented toward $u-1$, then $\tau^{-1}(u-1)>\tau^{-1}(v-1)$ and therefore $\ell(u-1,v-1)=2$ implying $\ell(u,v)=0$ and therefore it does
not matter that $\{u-1,v-1\}$ is not oriented from $u-1$ to $v-1$.

Finally assume $v<u$. Note that $v>0$ as otherwise the edge $\{u,v\}$ could not be directed toward $v$. If $\{u-1,v-1\}$ is not oriented from $u-1$ to
$v-1$, then the arc $\{u,u+1,\dots,v-1\}$ is either $A_{\tau^{-1}(u)}$ or
$A_{\tau^{-1}(v-1)}$. In the former case our rule implies that $\ell(u,v)=0$,
in the latter case it implies $\ell(u-1,v-1)=2$ and thus again $\ell(u,v)=0$
in which case we have no problem.
\medskip

\par\noindent
This proves that our rule gives a transitive $\sigma_n$-orientation and completes the
proof of the theorem.
\hfill$\Box$

\bigskip
\par\noindent
Note that the linear order obtained on the vertex set of $K_n$ by the
orientation in the proof above has some special properties. To formulate them
we introduce the following notions.

\begin{defi}
Let $\tau(1)\dots\tau(n)$ be an ordering of the numbers $0,1,\dots,n-1$. We say that $j\in\{0,1,\dots,n-1\}$
is a {\em local minimum} in this order if $j$ precedes both $j-1$ and $j+1$
(addition is meant modulo $n$), that is
$\tau^{-1}(j)<\tau^{-1}(j-1)$ and $\tau^{-1}(j)<\tau^{-1}(j+1)$.
Similarly,  $j\in\{0,1,\dots,n-1\}$
is a {\em local maximum} if $j$ is preceded by both $j-1$ and $j+1$, that is
$\tau^{-1}(j)>\tau^{-1}(j-1)$ and $\tau^{-1}(j)>\tau^{-1}(j+1)$.

We call $\tau$ {\em bitonic} if it has a unique local minimum and a unique local maximum.  We call a transitive $\sigma_n$-orientation of a $\sigma_n$-partition {\em standard} if it is consistent with a bitonic ordering of the vertices.
\end{defi}

\medskip
\par\noindent
\begin{prop}\label{felez}
The transitive $\sigma_n$-orientations given in the proof of Theorem~\ref{standard} are standard. The unique local minimum of the corresponding ordering is at $0$, while the unique local maximum is at $n/2$.
\end{prop}

\proof
The fact that $\tau$ given in the proof is bitonic and thus the transitive $\sigma_n$-orientations are standard follows immediately from the construction. It is also clear that $\tau(1)=0$ is the local minimum and we need only prove that the local maximum (that is $\tau(n)$) is $n/2$.

Recall from the proof of Theorem~\ref{standard} that $A_i=\{\tau(1),\dots,\tau(i)\}$ is a consecutive interval in the cycle formed by the $i$ vertices from $j_i+1$ to $j_i+i$ for all $1\le i\le n$. The sequence $(A_i)_{i=1}^n$ starts at $A_1=\{0\}$ and we obtain $A_{i+1}$ from $A_i$ by extending $A_i$ with $\tau(i+1)$ at one end of this interval. The label $\ell(j_i,j_i+i+1)$ (which is never $1$) determines which end we place $\tau(i+1)$, namely if the label is $0$ we chose one end, while if it is $2$ we chose the other end. For $1\le i\le n/2-1$ we have $\ell(j_i,j_i+i+1)\equiv a_{i+1}-a_i-1$ by Equation~(\ref{kicsi}) and $\ell(j_{n-i-1},j_{n-i-1}+n-i)\equiv a_i-a_{i+1}$ by  Equation~(\ref{nagy}). This means that  $\ell(j_i,j_i+i+1)+\ell(j_{n-i-1},j_{n-i-1}+n-i)\equiv2$, so one of these labels must be $0$ and the other $2$ and thus
we extend the interval $A_i$ on one end to get $A_{i+1}$ while we extend $A_{n-i-1}$ on the other end to obtain $A_{n-i}$. This partitions the $n-2$ extension steps bringing $A_1$ to $A_{n-1}$ into $n/2-1$ pairs and shows that we use $n/2-1$ extensions of the interval at either side. Thus $n/2$ (the vertex in distance $n/2$ from $A_1=\{0\}$ in either direction) must be the only vertex outside $A_{n-1}$ and therefore we have $\tau(n)=n/2$ as stated.
\hfill$\Box$
\medskip

The following proposition is a sort of converse of the previous one.

\begin{prop}\label{converse-standard}
If there exists a standard transitive $\sigma_n$-orientation for a $\sigma_n$-partition, then the conditions of Theorem~\ref{standard} are satisfied, namely the number $n$ of vertices is even and for the defining sequence $a_1\dots a_{n/2}$ of the $\sigma_n$-partition either $a_{j+1}=a_j$ or $a_{j+1}\equiv a_j+1 \pmod 3$ holds for each $1\le j\le n/2-1$. If there exists a standard transitive $\sigma_n$-orientation for a $\sigma_n$-partition, then it is unique up to shifts of the automorphism $\sigma_n^3$.
\end{prop}

\proof
Consider a standard transitive $\sigma_n$-orientation $T$ for a $\sigma_n$-partition. Without loss of generality we may assume that the defining sequence of the $\sigma_n$-partition starts with $a_1=0$. Consider the bitonic ordering $\tau$ of the vertices consistent with $T$. The first element $a=\tau(1)$ is clearly the local minimum. All edges are directed away from $a$, so $\sigma$ reverses the orientation of the edge $\{a-1,a\}$, as it brings it to $\{a,a+1\}$. Therefore we have $\ell(a,a+1)=0$ and thus $a\equiv0\pmod3$. We may and will assume $a=0$ as this can be achieved with a shift of a suitable power of the automorphism  $\sigma_n^3$.

As $\tau$ is bitonic, the set $A_i=\{\tau(1),\dots,\tau(i)\}$ must be an interval $\{j_i+1,\dots,j_i+i\}$ along the cycle formed by the vertices. Clearly, $\tau(i+1)$ is either $j_i$ or $j_i+i+1$. We show the uniqueness of $\tau$ by observing that the value of $\tau(i+1)$ depends on the label of the edge $e=\{j_i,j_i+i+1\}$ exactly as in the construction in the proof of Theorem~\ref{standard}, namely $\tau(i+1)=j_i$ if this label is $0$ and $\tau(i+1)=j_i+i+1$ if the label is $2$ and the label of $e$ cannot be $1$. Indeed, if $\tau(i+1)=j_i+i+1$, then $\sigma_n$ reverses the orientation of the edge $e$, so its label is $2$, but if $\tau(i+1)=j_i$, then $\sigma_n^{-1}$ reverses the orientation of $e$, so its label must be $0$.

From the rule established above we can derive Equation~(\ref{bar}) just as in the proof of Theorem~\ref{standard}. In case $n$ is odd we can apply this formula to $i=(n-1)/2$ and using also Equation~(\ref{shift}) we obtain $\ell(j_i,j_i+i+1)\equiv\ell(0,i+1)-\ell(0,i)-1=(\ell(0,i)+i+1)-\ell(0,i)-1\equiv i\equiv1$ contradicting that we saw that the label of the edge $\{j_i,j_i+i+1\}$ cannot be $1$. This proves that $n$ is even.

From Equation~(\ref{bar}) we can deduce  Equation~(\ref{kicsi}) for $1\le i\le n/2-1$. Using that the label of the edge $\{j_i,j_i+i+1\}$ cannot be $1$ this implies that $a_{i+1}=a_i$ or $a_{i+1}\equiv a_i+1$ finishing the proof of the proposition.
\hfill$\Box$

\subsection{Several parts}

In this subsection we extend the earlier results to partitions to more than three parts. The extensions are straightforward, the proofs carry over almost verbatim. We just list the results here and show how they apply to decompositions into Hamiltonian paths and cycles. Note that in the discussion below we also allow the case when the number of parts is $k=2$.

\begin{defi}
Let $k$ and $n$ be integers larger than $1$ and (as before) let $\sigma_n$ be the cyclic permutation on the set $[n]=\{0,1,\dots,n-1\}$ bringing $i$ to $i+1$. Recall that elements of $[n]$ are understood modulo $n$. Let us call a partition of the edge set of the complete graph on the vertex set $[n]$ into the subgraphs $F_0,\dots,F_{k-1}$ a {\em $\sigma_n$-$k$-partition} if $\sigma_n(F_d)=F_{d+1}$ for $0\le d<k-1$. This implies that $n$ is divisible by $k$, further if $k$ is even, then $n$ is divisible by $2k$. It also implies that $\sigma_n(F_{k-1})=F_0$ and that $\sigma_n^k$ is an automorphism of all the graphs $F_d$. We say that the {\em label} of an edge $e$ of $F_d$ is $\ell(e)=d$.
The defining sequence of this partition is the sequence $a_1,\ldots,a_{\lfloor n/2\rfloor}$, where $a_i=\ell(\{0,i\})$. The defining sequence uniquely determines the $\sigma_n$-$k$-partition, namely $F_d$ consists of the edges
$\{b,b+i\}$ for which $b\equiv d-a_i\pmod k$. This is indeed a $\sigma_n$-$k$-partition for any sequence $a_1,\dots,a_{\lfloor n/2\rfloor}$ over the letters $0,1,\dots,k-1$ if the divisibility conditions are satisfied by $n$ and $k$.

We call a transitive orientation $T$ of $K_n$ a {\em transitive $\sigma_n$-orientation} of a $\sigma_n$-$k$-partition if $\sigma_n(\vec{F_d})=\vec{F}_{d+1}$ for $0\le d<k-1$, where $\vec F_d$ is the orientation of $F_d$ obtained as a subgraph of $T$. A transitive orientation $T$ satisfies this condition if and only if $\sigma_n$ reverses the orientation of no edges of $F_d$ with $d\ne k-1$. We call a transitive $\sigma_n$-orientation {\em standard} if it is consistent with a bitonic ordering of the vertices.
\end{defi}

The common generalization of Theorems~\ref{thm:selfc} and \ref{standard} is as
follows.

\begin{thm}\label{genstandard}
There exists a standard orientation for a $\sigma_n$-$k$-partition if and only if $n$ is divisible by $2k$ and the defining sequence $a_1a_2\dots a_{n/2}$ of the partition satisfies that for every $1\le j\le n/2-1$ either $a_{j+1}=a_j$ or $a_{j+1}\equiv a_j+1\pmod k$. If a standard orientation exists, then it is unique up to shifts of $\sigma^k$ and the corresponding bitonic ordering of the vertices satisfies that the first and last elements in the ordering (the unique local minimum and maximum) differ by $n/2$.
\end{thm}
\medskip

\par\noindent
\begin{remark}
It is straightforward to see that Theorem~\ref{genstandard} generalizes
Theorem~\ref{standard}. Note that Theorem~\ref{thm:selfc} is also implied by
Theorem~\ref{genstandard} since in case of two self-complementary graphs, that
is, when $k=2$, the condition $a_{j+1}\equiv a_j\ \ {\rm or}\ \  a_j+1\pmod k$
is automatically satisfied as $a_{j+1}$ cannot take any value other than $0$
or $1$.
\end{remark}
\medskip

We formulate the following simple generalization of the trivial observation that the (only) $\sigma_3$-partition has no transitive $\sigma_3$-orientation.

\begin{prop}\label{trivi}
If there exists a transitive $\sigma_n$-orientation for a $\sigma_n$-$k$-partition, then $n\ge2k$ and $n\ne3k$. If $n=2k$, then any transitive $\sigma_n$-orientation of a $\sigma_n$-$k$-partition is standard.
\end{prop}

\proof
From the existence of this partition we know that $k$ divides $n$. Note that $\sigma_n$ reverses the orientation of an even number of the $n$ edges $\{i,i+1\}$ in any orientation $T$ of $K_n$. But if $k=n$, then only one of these edges has label $k-1$, so if $T$ is a transitive $\sigma_n$-orientation, then the orientation of none of these edges are reversed. But then they form a directed cycle, so $T$ is not transitive.

If $n=2k$ or $3k$, then two or three of these $n$ edges have label $k-1$, so exactly two of them have to be reversed by $\sigma_n$ to avoid the directed cycle. This makes the transitive $\sigma_n$-orientation standard. But we know from Theorem~\ref{genstandard} that if a $\sigma_n$-$k$-partition admits a standard orientation, then $n$ is divisible by $2k$, so $n=3k$ is not an option.
\hfill$\Box$
\medskip

Complete graphs of odd order can be decomposed into Hamiltonian cycles. This
is a classical result as mentioned in Adrian Bondy's chapter \cite{Bondy} of
the {\em Handbook of Combinatorics}. (It is added there that ``one such
construction, due to a Monsieur Walecki, is described in the book by Lucas
(1891, pp. 161--164)'', cf. \cite{Lucas}.)

This result extends to the decomposition of odd order transitive tournaments to identically oriented Hamiltonian cycles. This decomposition can be found directly, but it can also be arrived at by applying our result to a certain decomposition of the non-oriented complete graph as shown below.

\begin{cor}\label{Hpaths}
If $n$ is even, then we can decompose $T_n$ into $\frac{n}{2}$ alternatingly
oriented Hamiltonian paths.
\end{cor}

\proof
Let $F_0$ be the path defined as
$0,1,n-1,2,n-2,\dots,i,n-i,(i+1),\dots,(\frac{n}{2}+1),\frac{n}{2}$. It is easy to see that the graphs
$F_d:=\sigma_n^d(F_0)$ for $d=0,1,\dots,n/2-1$
partition the edge set of $K_n$ into $\frac{n}{2}$ Hamiltonian paths. This is a $\sigma_n$-$n/2$-partition.

One can also readily check that $\{0,j\}\in E(F_{\lfloor j/2\rfloor})$ holds for every
$j\in\{1,\dots,n-1\}$. This means that the defining sequence of this
$\sigma$-$n/2$-partition is $a_1,\dots,a_{n/2}$ with $a_i=\lfloor
i/2\rfloor$. By Theorem~\ref{genstandard} there is a standard orientation for
this partition that is unique up to a shift with $\sigma_n^{n/2}$. It is not
hard to check that the standard orientation orients the Hamiltonian paths in
this partition alternatingly, that is, each vertex will be a source or a sink of each Hamiltionian path. This construction is illustrated for
$n=6$ by the oriented path in the very central picture of Figure~1.
\hfill$\Box$

\begin{cor}\label{Hcycles}
If $n$ is odd, then we can decompose $T_n$ into $\frac{n-1}{2}$ isomorphically
oriented Hamiltonian cycles.
\end{cor}

\proof
Consider the decomposition of $T_{n-1}$ given in Corollary~\ref{Hpaths} on
vertices labeled by $[n-1]$. Add the extra vertex $v$ and
connect it to the two endpoints of each of the Hamiltonian paths thus extending
them to Hamiltonian cycles. Orient all edges incident to $v$ away from $v$. These isomorphically oriented Hamiltonian cycles decompose the transitive tournament on $n$ vertices.
\hfill$\Box$.

\begin{remark}
The orientation of the Hamiltonian cycles in our decomposition in
Corollary~\ref{Hcycles} is such that all but one of the vertices is either a source or a sink. This kind of orientation of odd cycles is called
alternating in \cite{KPS}, where it is shown that these oriented versions of
odd cycles have maximal Sperner capacity. In the special case of $n=5$ this
orientation already appeared in \cite{GGKS}, where it was observed that its
Sperner capacity is $\sqrt{5}$, that is, it achieves the Shannon capacity of the
underlying undirected graph which is $C_5$, whose capacity was determined in
the celebrated paper by Lov\'asz \cite{LL}. This observation was the starting
point of our investigations in \cite{SaSi}.

In Section~\ref{sect:smallex} we explained how to find transitive $\sigma$-orientations for $\sigma$-partitions if the permutation $\sigma$ is not cyclic. First we solve the restriction of the problem for the domain of each cycle in the cycle decomposition of $\sigma$, then extend the obtained orientations by orienting the edges between distinct cycles consistently with a linear ordering of these cycles. We used the same strategy here for the two cycles of the permutation $\sigma$ that has $v$ as a fixed point and acts on $[n-1]$ as $\sigma_{n-1}$. The decomposition has $v$ as the domain of a trivial cycle and $[n-1]$ as the domain of $\sigma_{n-1}$. Had we defined these notions for arbitrary permutations (not just for cyclic ones), we could call the decomposition of the complete graph in Corollary~\ref{Hcycles} a $\sigma$-$(n-1)/2$-partition of $K_n$ and the orientations of the Hamiltonian cycles would form a transitive $\sigma$-orientation of this partition.
\end{remark}

\section{Non-standard orientations}

The results in the previous section may make one hope that the conditions of Theorem~\ref{standard} are not just sufficient but also necessary for a $\sigma_n$-partition to have transitive $\sigma_n$-orientations. This is not the case.
In this section we give some sufficient conditions that go beyond the ones in Theorems~\ref{standard} and \ref{genstandard} .

Our construction takes a transitive $\sigma_m$-orientation for a
$\sigma_m$-$k$-partition and uses that to get transitive
$\sigma_n$-orientations of related $\sigma_n$-$k$-partitions. Here $n$ is a multiple of $m$. Even for some standard $\sigma_m$-$k$-orientations, the resulting transitive $\sigma_n$-orientations are not always standard and in some cases the $\sigma_n$-$k$-partitions do not have standard orientations at all.


\begin{defi}\label{defi:blowup}
Let $k>1$ and $m$ be such that $\sigma_m$-$k$-partitions exist (that is, $m$ is a multiple of $k$ and if $k$ is even $m$ is also a multiple of $2k$). Let $n$ be a multiple of $m$. We call the $\sigma_n$-$k$-partition $P$ with defining sequence $a_1,\dots,a_{\lfloor n/2\rfloor}$ a {\em blow-up} of the $\sigma_m$-$k$-partition $Q$ with defining sequence $b_1,\ldots,b_{\lfloor m/2\rfloor}$ if $a_i=b_{(i\bmod m)}$ whenever $1\le i\le\lfloor n/2\rfloor$ and $i$ is not divisible by $m$. For this to make sense even if $i\bmod m>\lfloor m/2\rfloor$ we extend the sequence $b_i$ by setting $b_i=(b_{m-i}+i)\bmod k$ for $\lfloor m/2\rfloor<i<m$. This makes $b_i$ the label of the edge $\{0,i\}$ in the $\sigma_m$-$k$-partition $Q$ for any $i$. Note that we have no requirement for the value of $a_i$ if $i$ is divisible by $m$, so $P$ is not determined by $Q$ and $n$.
\end{defi}

Notice that if the $\sigma_m$-$k$-partition $Q$ has a standard orientation, then its defining sequence satisfies the requirements of Theorem~\ref{genstandard} and therefore the defining sequence $a_1,\dots,a_{n/2}$ of its blow-up $P$ also satisfies $a_{i+1}=a_i$ or $a_{i+1}\equiv a_i+1\pmod k$ whenever neither $i$
nor $i+1$ is divisible by $m$. However, by the free choice of the value $a_i$
whenever $i$ is divisible by $m$, this property need not hold for
the other indices, thus $P$ may violate the
conditions of Theorem~\ref{genstandard}. Nevertheless, as we will show, $P$ admits a transitive $\sigma_n$-orientation in this case. We call the transitive orientation constructed in the next theorem the {\em blow-up} of the transitive $\sigma_m$-orientation of $Q$.

\begin{thm}\label{thm:blowup}
If the $\sigma_n$-$k$-partition $P$ is a blow-up of the $\sigma_m$-$k$-partition $Q$ and $Q$ admits a transitive $\sigma_m$-orientation, then $P$ admits a transitive $\sigma_n$-orientation.
\end{thm}
\medskip

\par\noindent
\proof
Let $T$ be the transitive $\sigma_m$-orientation of $Q$ we assumed to exist and let $\tau$ be the ordering of the vertex set $[m]$ that $T$ is consistent with. The theorem claims that $P$ has transitive $\sigma_n$-orientation $T'$. We construct $T'$ by finding the ordering $\tau'$ on the vertex set $[n]$ that $T'$ is consistent with. The only requirement we have to satisfy is that all edges of $T'$ whose orientation $\sigma_n$ reverses should have label $k-1$.

Let us set $d=n/m$ and for $i\in[m]$, let $H_i=\{jm+i\mid j\in[d]\}$. Each of these $m$ sets has $d$ elements and together they partition $[n]$. The ordering $\tau'$ starts with the elements of $H_{\tau(1)}$ followed by the elements of $H_{\tau(2)},\dots,H_{\tau(m-1)}$. The order within the sets $H_i$ will be specified later. We call an edge of $K_n$ an {\em outer edge} if it connects vertices from distinct sets $H_i$ and $H_j$, otherwise it is an {\em inner edge}. The orientation of the outer edges in $T'$ are not influenced by the order within the sets $H_i$. namely for $a\in H_i$ and $b\in H_j$ ($i\ne j$) the orientation of the edge $\{a,b\}$ in $T'$ is determined by the orientation of the of edge $\{i,j\}$ in $T$. This means that $\sigma_n$ reverses the orientation of the edge $\{a,b\}$ in $T'$ if and only if $\sigma_m$ reverses the orientation of the edge $\{i,j\}$ in $T$. Our definition of a blow-up ensures that the label of the edge $\{a,b\}$ for the partition $P$ is the same as the label of the edge $\{i,j\}$ for the partition $Q$. Therefore $\sigma_n$ reverses the orientation of outer edges in $T'$ only if their label is $k-1$.

To finish the proof we have to specify the ordering $\tau'$ within the sets $H_i$ in such a way that the same can be said about inner edges of $T'$: $\sigma_n$ reserves the orientation of them only if their label is $k-1$.

For $i\in[m]$ let $\tau_i(1),\tau_i(2),\dots,\tau_i(d)$ be an ordering of $[d]$ to be specified later and let us say that $\tau'$ orders the elements of $H_i$ in the following order: $\tau_i(1)m+i,\tau_i(2)m+i,\dots,\tau_i(d)m+i$.

Consider an inner edge $e=\{a,b\}$, with $a,b\in H_i$, $i\in[m]$. We have $a=Am+i$ and $b=Bm+i$ for some $A,B\in[d]$ and the orientation of $e$ in $T'$ is determined by the order of $A$ and $B$ in $\tau_i$. In case $i\ne m-1$ we have $\sigma_n(e)=\{a+1,b+1\}$, $a+1,b+1\in H_{i+1}$ and the orientation of $\sigma_n(e)$ is determined by the order of $A$ and $B$ in $\tau_{i+1}$. Thus, $\sigma_n$ reverses the orientation of $e$ in $T'$ if and only if $A$ and $B$ are in different order in $\tau_i$ and in $\tau_{i+1}$. We must make sure that this only happens if the label of $e$ is $k-1$.

The situation is a bit different for inner edges $e=\{a,b\}$ with $a,b\in H_{m-1}$. We still have $\sigma_n(e)=\{a+1,b+1\}$ (recall that these vertices are understood modulo $n$). But now $a+1,b+1\in H_0$ and $a+1=(A+1)m$, $b+1=(B+1)m$. Thus, $\sigma_n$ reverses the orientation of $e$ if and only if $A$ and $B$ are in different order in $\tau_{m-1}$, then $A+1$ and $B+1$ are in $\tau_0$. Here both $A+1$ and $B+1$ are understand modulo $d$. We will make sure that this only happens if the label of $e$ is $k-1$.

We will now specify the orderings $\tau_i$. As explained above, this finishes the proof if we can show the following two properties:
\begin{description}
\item{(a)} If for some $0\le i<m-1$ and $A\ne B\in[d]$ the order of $A$ and $B$ is different in $\tau_i$ and $\tau_{i+1}$, then we have $\ell(Am+i,Bm+i)=k-1$.
\item{(b)} If $A$ precedes $B$ in $\tau_{m-1}$ but $(B+1)\bmod d$ precedes $(A+1)\bmod d$ in $\tau_0$, then we have $\ell(Am+m-1,Bm+m-1)=k-1$.
\end{description}

We set $\tau_0$ to be the order: $0,1,d-1,2,d-2,\dots,\lceil d/2\rceil$. Note that the pair of elements $j$ and $d-j$ are consecutive for all $1\le j<d/2$.  For $0\le i<k$ we obtain $\tau_{i+1}$ from $\tau_i$ by swapping the order of $j$ and $d-j$ for each $j$ satisfying $\ell(jm+i,(d-j)m+i)=k-1$. In this way we maintain that the $j$ and $d-j$ are consecutive in $\tau_i$ for all $0\le i\le k$ and $1\le j<d/2$. Also, with this rule condition (a) is satisfied for $0\le i<k$. For any $i$ and $j$ as above we have $\ell(jm+i,(d-j)m+i)\equiv\ell(jm,(d-j)m)+i\pmod k$. Therefore, for any such $j$, the label will be $k-1$ for exactly one of the indices $0\le i<k$ and thus $\tau_k$ will have all the pairs $(j,d-j)$ swapped. Namely, $\tau_k$ is the order $0,d-1,1,d-2,2,d-3,\ldots,\lfloor d/2\rfloor$.

Proposition~\ref{trivi} tells us that $m\ge2k$. Let us assume for now that $m>2k$. We will come back to the case $m=2k$ later. Observe that $j$ and $d-j-1$ are consecutive in $\tau_k$ for any $0\le j\le d/2-1$. For $k\le i<2k$ we obtain $\tau_{i+1}$ from $\tau_i$ by swapping the order of $j$ and $d-j-1$ for each $j$ satisfying $\ell(jm+i,(d-j-1)m+i)=k-1$. In this way we maintain that the vertices $j$ and $d-j-1$ are consecutive in $\tau_i$ for all $k\le i\le2k$ and $0\le j\le d/2-1$. Also, this rule makes condition (a) satisfied for $k\le i<2k$. Just as before, for any $j$ as above the label condition is satisfied for exactly one index $k\le i<2k$ and thus $\tau_{2k}$ will have all the pairs $(j,d-j-1)$ swapped. Namely, $\tau_{2k}$ is the order $d-1,0,d-2,1,d-3,2,\ldots\lceil d/2\rceil-1$.

We set $\tau_i=\tau_{2k}$ for $2k<i<m$. This makes condition (a) hold vacuously for $2k\le i\le m-1$ as $\tau_i=\tau_{i+1}$. Condition (b) is also satisfied vacuously, since $A$ precedes $B$ in $\tau_{m-1}=\tau_{2k}$ if and only if $(A+1)\bmod d$ precedes $(B+1)\bmod d$ in $\tau_0$. This is so because $\tau_0$ can be obtained from $\tau_{2k}$ by replacing each element $j$ by $(j+1)\bmod d$. This finishes the proof of the theorem in the case $m>2k$.

It remains to consider the case $m=2k$. We define the orders $\tau_i$ for $k<i\le2k$ the same way as above. We do not use $\tau_{2k}$ in the definition of $\tau'$, but we will use it in our argument below.

Condition (a) is satisfied as before. But now condition (b) is not vacuous. We still have that $A$ precedes $B$ in $\tau_{2k}$ if and only if $(A+1)\bmod d$ precedes $(B+1)\bmod d$ in $\tau_0$, so if $A$ precedes $B$ in $\tau_{2k-1}$ but $(B+1)\bmod d$ precedes $(A+1)\bmod d$ in $\tau_0$ as called for in condition (b), then $A$ and $B$ appear in different order in $\tau_{2k-1}$ and $\tau_{2k}$ and therefore $\ell(Am+m-1,Bm+m-1)=k-1$. This makes condition (b) satisfied and finishes the proof of the theorem.
\hfill$\Box$

To illustrate Theorem~\ref{thm:blowup} we show a $\sigma_n$-partition that admits a transitive $\sigma_n$-orientation but does not admit a standard one. We take $n=12$ and consider the $\sigma$-partition $P$ with defining sequence $000121$. Recall that this is a partition of the edge set of $K_{12}$ to three isomorphic subgraphs and it has no standard orientation by Proposition~\ref{converse-standard}. But it is a blow-up of the $\sigma_6$-partition $Q$ with defining sequence $000$. $Q$ has a transitive $\sigma_6$-orientation, even a standard one by Theorem~\ref{standard}. The first part of $Q$ and its orientation is depicted in the first illustration on Figure~1. By Theorem~\ref{thm:blowup} $P$ has a transitive $\sigma_{12}$-orientation. For example orienting the edges consistent with the following ordering of the vertices gives a transitive $\sigma_{12}$-orientation: $0,6,1,7,2,8,11,5,4,10,3,9$.

\section{Complications}

After seeing Theorems~\ref{standard} and \ref{thm:blowup} one might be curious
to know whether they describe all transitive $\sigma_n$-orientations of a
$\sigma_n$-partition, namely if all such orientations are blow-ups of standard
orientations. This is, however, not the case, moreover, there exist
$\sigma_n$-partitions that neither admit a standard orientation nor are
blow-ups of $\sigma_m$-partitions for some $m<n$, yet they do admit a
transitive $\sigma_n$-orientation. We found such examples by computer and do
not see a general pattern that would still suggest a complete characterization.
(We note that our examples in this section all concern the simplest possible case $k=3$ again.)

The $\sigma_{24}$-partition with the defining sequence
$$0,0,0,1,2,0,0,0,1,1,2,1$$
is such an example. The orientation consistent with the following order of the vertices is a transitive $\sigma_{24}$-orientation for this partition:

$$0,1,2,23,22,21,3,9,4,10,20,5,11,8,7,19,6,18,12,13,14,17,16,15.$$

(To check that this defines a transitive $\sigma_{24}$-orientation one has only to verify that the edges whose orientation $\sigma_{24}$ reverses have label $2$. For example, the edge $e=\{4,8\}$ is oriented towards $8$ as $8$ appears later in this sequence than $4$. But $\sigma_{24}(e)=\{5,9\}$ is oriented towards $5$ as $5$ appears later than $9$. So $\sigma_{24}$ reverses the orientation of $e$. Now $\ell(4,8)\equiv\ell(0,4)+4\pmod3$ and $\ell(0,4)$ is the fourth number of the defining sequence, namely $1$, so $\ell(4,8)=2$ as required.)

\section{Necessary conditions}

So far we have seen sufficient conditions for $\sigma_n$-partitions (or $\sigma_n$-$k$-partitions) to posses a transitive $\sigma_n$-orientation. While a complete characterization seems elusive it makes sense to look also for non-trivial necessary conditions.
Here we give a simple such condition.

\begin{thm}\label{prop:kkkz} Let $P$ be a $\sigma_n$-$k$-partition with defining sequence $a_1,\ldots,a_{\lfloor n/2\rfloor}$. Assume $a_1=0$ and let $i$ be an index with $1\le i<\lfloor n/2\rfloor$ such that $a_j\ne k-1$ for $1\le j\le i$. If $P$ has a transitive $\sigma_n$-orientation, then $a_{i+1}\le a_i+1$.
\end{thm}

\proof
Fix a transitive $\sigma_n$-orientation $T$ of $P$. When referring to the orientation of edges we speak about the orientation in $T$.

We call a vertex $m\in[n]$ a {\em leader} if it is divisible by $k$. Notice that if $m$ is a leader, then the label $\ell(m,m+j)$ is $a_j$ for all $1\le j\le\lfloor n/2\rfloor$. (Recall that vertices are always understood modulo $n$.) We call a leader $m$ an {\em in-leader} if the edge $\{m,m+1\}$ is oriented towards $m$, otherwise it is an {\em out-leader}.

We claim that if $m$ is an in-leader, then all the edges $\{m,m+j\}$ for $1\le j\le i+1$ are oriented toward $m$, while if $m$ is an out-leader, then all these edges are oriented away from $m$. By symmetry, it is enough to prove one of the statements. We prove the latter one by induction on $j$. The statement of the claim is assumed for $j=1$. So let $1\le j\le i$ and assume $\{m,m+j\}$ is oriented away from $m$. The label of this edge is $a_j\ne k-1$, so $\sigma_n$ does not reverse its orientation. Therefore the edge $\{m+1,m+j+1\}$ is oriented towards $m+j+1$. As $m$ is an out-vertex the edge $\{m,m+1\}$ is oriented toward $m+1$. To get a transitive orientation $\{m,m+j+1\}$ must therefore be oriented away from $m$ finishing the inductive proof.

Consider the cycle formed by the edges $\{m,m+1\}$ for all $m\in[n]$. If all leader vertices were in-leaders or all of them were out-leaders, then this cycle would be a directed cycle contradicting the transitivity of $T$. So we must have at least one in-leader and also at least one out-leader. We can therefore fix an out-leader $m$ such that the very next leader vertex, namely $m+k$ is an in-leader.

By the claim above (and since $m$ is an out-leader) the edge $\{m,m+i\}$ is oriented towards $m+i$. The label of this edge is $a_i$, so one can apply the permutation $\sigma_n$ to this edge $k-1-a_i$ times without it reversing the orientation. Thus, the edge $\{m+k-1-a_i,m+i+k-1-a_i\}$ is oriented toward $m+i+k-1-a_i$. Recall that $a_i\ne k-1$, so (as $m$ is an out-leader) the edge $\{m+k-2-a_i,m+k-1-a_i\}$ is oriented toward $m+k-1-a_i$. The transitivity of $T$ therefore implies that the edge $e=\{m+k-2-a_i,m+i+k-1-a_i\}$ is oriented toward $m+i+k-1-a_i$. We have $\sigma_n^{a_i+2}(e)=\{m+k,m+k+i+1\}$ and this edge is oriented toward the in-leader $m+k$ by the claim above. We see that $\sigma_n^{a_i+2}$ reverses the orientation of $e$, therefore the label $a_{i+1}$ of $\sigma_n^{a_i+2}(e)$ must be less than $a_i+2$.
\hfill$\Box$
\medskip

\begin{defi}
Let $a_1,\ldots,a_{\lfloor n/2\rfloor}$ be the defining sequence of a $\sigma_n$-$k$-partition $P$. We say that the sequence {\em halts} at the index $i$ ($1\le i<\lfloor n/2\rfloor$) if $a_{i+1}=a_i$. We say that it {\em steps} at the index $i$ if $a_{i+1}\equiv a_i+1\pmod k$. We say that it {\em jumps} at the index $i$ if it neither halts nor steps there.

 We call the $\sigma_n$-$k$-partition with defining sequence $b_1,\ldots,b_{\lfloor n/2\rfloor}$ the {\em dual} of $P$ if $b_i\equiv i-1-a_i\pmod k$ for all $i$.
\end{defi}

Note that if a the defining sequence of a $\sigma_n$-$k$-partition halts at an index $i$, then the defining sequence of its dual steps there and vice versa. The defining sequences of a $\sigma_n$-$k$-partition and its dual jump at the same indices.

We can rephrase Theorem~\ref{genstandard} as follows: A $\sigma_n$-$k$-partition has a standard orientation if and only if its defining sequence does not jump at all.

\begin{thm}\label{dual}
  If a $\sigma_n$-$k$-partition admits a transitive $\sigma_n$-orientation, then so does its dual.

  If a $\sigma_n$-$k$-partition $P$ admits a transitive $\sigma_n$-orientation and the defining sequence of $P$ jumps at an index $i$, then there is an index $j<i$ where it halts and at least $k-1$ distinct indices $j'<i$ where it steps, or the other way around: there is an index $j<i$ where it steps and at least $k-1$ distinct indices $j'<i$ where it halts. In particular, there is no jump at indices $i\le k$.
\end{thm}

\proof
Let $P$ be the $\sigma_n$-$k$-partition, with parts $F_0,\dots,F_{k-1}$. The dual $Q$ of $P$ can be obtained by first relabeling the vertices, namely switching the label $v$ and $n-v$ for $1\le v<k$, and then considering the parts in reverse order, namely as $F_{k-1},F_{k-2},\dots,F_0$. If a transitive orientation $T$ of $K_n$ is a transitive $\sigma_n$-orientation of $P$, then $T$ (after the relabeling) is also a transitive $\sigma_n$-orientation of $Q$. This proves the first claim in the theorem.

To verify the second claim we assume without loss of generality that the defining sequence $(a_j)$ of $P$ starts with $a_1=0$. Note that this makes the defining sequence $(b_j)$ of its dual $Q$ start with $b_1=0$. Assume that the first jump in $(a_j)$ is at the index $i$, so at indices $1\le j<i$ the sequence halts or steps.
By Theorem~\ref{prop:kkkz} we have $a_{i+1}\le a_i+1$ unless the sequence steps for at least $k-1$ such indices. Further, it has to step for at least one such index, as otherwise we would have $a_j=0$ for all $j\le i$ and $a_{i+1}\le a_i+1=1$ contradicting our assumption that the sequence jumps at $i$.
Note that $i$ is also the first index where the sequence $(b_j)$ jumps, so we similarly have that this sequence steps for at least one index $j<i$ and we have $b_{i+1}\le b_i+1$ unless it steps at $k-1$ or more such indices. To finish the proof it is enough to note that $(a_j)$ steps where $(b_j)$ halts and vice versa and if both step at fewer than $k-1$ indices $j<i$, then $a_i\le a_{i+1}\le a_i+1$ contradicting our assumption that the sequence $(a_j)$ jumps at $i$.
\hfill$\Box$
\medskip

\noindent
We conclude the paper with two conjectures.
All $\sigma_n$-partitions for which we found a transitive $\sigma_n$-orientation were for even values of $n$. We tried to find such examples with odd $n$ but failed even with a computer. This suggests the following.

\begin{conj}
No $\sigma_n$-partition with $n$ odd has a transitive $\sigma_n$-orientation.
\end{conj}

Although we did not do any computer search for partitions with more than three parts, we still venture the following stronger conjecture:

\begin{conj}
No $\sigma_n$-$k$-partition with $k>1$ and $n$ odd has a transitive $\sigma_n$-orientation.
\end{conj}


\begin{thebibliography}{99}

\bibitem{Bondy} J.\ A. Bondy, Basic graph theory: paths and circuits, in: {\em
  Handbook of combinatorics, Vol. 1}, 3--110, Elsevier Sci. B. V.,
  Amsterdam, 1995.

\bibitem{GGKS} A. Galluccio, L. Gargano, J. K\"orner, G. Simonyi, Different
  capacities of a digraph, {\em Graphs Combin.}, {\bf 10} (1994), 105--121.

\bibitem{english} R.\ A.\ Gibbs, Self-complementary graphs, {\sl J. Combin.
Theory ser. B}, {\bf 16} (1974), 106--123.

\bibitem{GKMPW}
A. G\"orlich, R. Kalinowski, M. Meszka, M. Pilsniak, M. Wo\'zniak,
A note on decompositions of transitive tournaments.
{\em Discrete Math.} 307 (2007), no. 7-8, 896�904.

\bibitem{GKMPW2}
A. G\"orlich, R. Kalinowski, M. Meszka, M. Pilsniak, M. Wo\'zniak,
A note on decompositions of transitive tournaments. II.
{\em Australas. J. Combin.} 37 (2007), 57�66.

\bibitem{GP} A. G\"orlich, M. Pil\'sniak,
A note on an embedding problem in transitive tournaments,
{\em Discrete Math.} 310 (2010), no. 4, 681�686.

\bibitem{GPW}A. G\"orlich, M. Pilsniak, M. Wo\'zniak,
A note on a packing problem in transitive tournaments,
{\em Graphs Combin.} 22 (2006), no. 2, 233�239.

\bibitem{GyA}  A.\ Gy\'arf\'as, Transitive tournaments and self-complementary
  graphs. {\em J. Graph Theory}, {\bf 38} (2001), no. 2, 111--112.

\bibitem{KPS} J. K\"orner, C. Pilotto, G. Simonyi, Local chromatic number and
Sperner capacity, {\em J. Combin. Theory Ser. B}, {\bf 95} (2005), 101--117.

\bibitem{LL}  L.\ Lov\'asz, On the Shannon capacity of a graph, {\em IEEE Trans.
Inform. Theory}, {\bf 25} (1979), 1--7.

\bibitem{Lucas} \'E.\ Lucas, {\em R\'ecr\'eations Mathematiques}, Vol. 2,
  Gauthiers-Villars, Paris, 1891.

\bibitem{Ringel} G.\ Ringel, Selbstkomplement\"are Graphen, {\em Arch.
Math.}, {\bf 14} (1963), 354-358.

\bibitem{Sachs} H.\ Sachs, \"Uber selbstkomplement\"are Graphen, {\em
  Publ. Math. Debrecen}, {\bf 9} (1962), 270--288.

\bibitem{SaSi} A. Sali, G. Simonyi, Orientations of self-complementary graphs
  and the relation of Sperner and Shannon capacities, {\em Europ.
    J. Combin.}, {\bf 20} (1999), 93--99.









\end{thebibliography}
\end{document}